\documentclass[10pt,twoside]{article}
\usepackage{amsmath}
\usepackage{Latex-document}

\newtheorem {theorem}{Theorem}

\def\eps{\epsilon}
\def\ep{\epsilon}
\def \rar{\rightarrow}
\def\R{{\bf R}}
\def\T{{\bf T}}
\def\P{{\bf P}}
\def\E{{\bf E}}
\def\Z{\hbox{\bf Z}}
\def \TT{{\mathcal T}}
\def \CC{{\mathcal C}}
\def \C{{\mathcal C}}

\title{\bf  Brownian Intersections, Cover Times \vskip -2mm  and Thick Points via Trees \vskip 6mm}

\author{Yuval Peres\vspace*{-0.5cm}\thanks{Departments of Statistics
\& Mathematics, University of California, Berkeley CA, USA.
 E-mail: peres@stat.berkeley.edu}}
\date{\vspace{-8mm}}

\begin{document}

\maketitle

\thispagestyle{first} \setcounter{page}{73}

\begin{abstract}\vskip 3mm

There is a close connection between intersections of Brownian motion
paths and percolation on trees.
 Recently, ideas from probability on trees
were an important component of the multifractal analysis of Brownian
occupation measure, in joint work with A. Dembo, J. Rosen and
O. Zeitouni. As a consequence, we proved two conjectures about
simple random walk in two dimensions: The first, due to Erd\H{o}s
and Taylor (1960), involves the number of visits to the most
visited lattice site in the first $n$ steps of the walk. The
second, due to Aldous (1989), concerns the number of steps it
takes a simple random walk to cover all points of the $n$ by
$n$ lattice torus. The goal of the lecture is to relate how
methods from probability on trees can be applied to random walks
and Brownian motion in Euclidean space.

\vskip 4.5mm



\end{abstract}

\vskip 12mm

\section{Introduction} \label{section 1}\setzero
\vskip-5mm \hspace{5mm }

In \cite{in}, the author showed that long-range
intersection probabilities for random
walks, Brownian motion paths and Wiener sausages in Euclidean space, can be
estimated up to constant factors by
survival probabilities of percolation processes on trees.

More recently, several long-standing problems involving
cover times and ``thick points'' for random walks in
two dimensions were solved in joint works \cite{cover, planar} of A. Dembo,
J. Rosen, O. Zeitouni and the author. These solutions were motivated
by powerful analogies with corresponding problems on trees,
but these analogies were not discussed explicitly in the
research papers cited. The goal of the
present note is to describe the tree problems and solutions,
that correspond to the problems studied in \cite{cover,planar}.

The {\em cover time} for a random walk on a finite graph
is the number of steps it takes the random walk to visit all vertices.
The cover time has been studied intensively by probabilists,
combinatorialists, statistical physicists and computer scientists,
with a variety of motivations; see, e.g.,
\cite{bro, mihail, aldous2,hil,nemi}.
The problem of determining the expected cover time $\TT_n$ for
the $n$ by $n$ lattice torus $\Z_n^2$,
was posed by Wilf~\cite{wilf} and Aldous~\cite{aldous1}.
In \cite{cover} we proved the following conjecture of Aldous~\cite{aldous1}.
\begin{theorem}
\label{theo-2lattice}
If $\TT_n$ denotes the time it takes
for the simple random walk in $\Z_n^2$ to completely cover $\Z_n^2$, then
\begin{equation} \label{eq:covplane}
\lim_{n\to \infty} \frac{\TT_n}{(n\log n)^2}=
\frac{4}{\pi}
\hspace{.1in}\mbox{in probability.}
\end{equation}
\end{theorem}

The first step toward proving Theorem \ref{theo-2lattice},
was to find a sufficiently robust proof for the asymptotics
of the cover time of finite $b$-ary trees. These asymptotics were
originally determined by Aldous in \cite{aldous-t}, but
his elegant recursive method was quite sensitive
and did not adapt to the approximate tree structure
that can be found in Euclidean space.
Cover times on trees are discussed in the next section.

Turning to a different but related topic,
Erd\H{o}s and Taylor (1960) posed a problem about
simple random walks in $\Z^2$: \em
How many times does the walk revisit the most
frequently visited site in the first $n$ steps? \rm

\begin{theorem}[\cite{planar}]
\label{ETlattice}
Denote by
$T_n(x)$ the number of visits of planar simple random walk
to $x \in \Z^2$ by time $n$, and let $T_n^* :=\max_{x\in \Z^2} T_n(x)\,.$
Then
\begin{equation}
\label{ErdT1}
\lim_{n\to\infty} \frac{T_n^*}{(\log n)^2} =
 \frac{1}{ \pi }\, \quad a.s.\,.
\end{equation}
\end{theorem}

This was conjectured by Erd\H{o}s and Taylor~\cite[(3.11)]{ET}.
After D. Aldous heard one of us describing this result,
he pointed us to his cover time conjecture,
and this eventually led to Theorem \ref{theo-2lattice}.
Although the proofs of that theorem and of Theorem \ref{ETlattice}
differ in important technical points, they follow
the same basic pattern:

\begin{description}
\item{(i)} Formulate  a suitable tree-analog and find a ``robust'' proof.
\item{(ii)} Establish a Brownian version using excursion counts.
\item{(iii)} Deduce the lattice result via strong approximation
a-la \cite{KMT}.
\end{description}

\section{Cover times for trees} \label{section 2} \setzero
\vskip-5mm \hspace{5mm}

Let $\Gamma_k$ denote the balanced $b$-ary tree of height $k$,
which has \newline $n_k=(b^{k+1}-1)/(b-1)$ vertices, and $n_k-1$ edges.

\begin{theorem}[Aldous~\cite{aldous-t}]
\label{cover-bary}
Denote by $\C_k$ the time it takes
for simple random walk in $\Gamma_k$, started at the root, to cover $\Gamma_k$. Then
\begin{equation} \label{eq:covtree}
\lim_{k\to \infty}\frac{\E\C_k}{n_k k^2}=
2\log(b) \, .
\end{equation}
\end{theorem}

\noindent{\bf Remark}  The expected hitting time from one vertex
to another is bounded by the commute time, which equals the
effective resistance times twice the number of edges (see, e.g.,
\cite{Ald-F}). Therefore the expected hitting time between two
vertices in $\Gamma_k$ is at most $4kn_k$. From a general result
in \cite{aldous3}, it follows that also
\begin{equation} \label{eq:covtree2}
\lim_{k\to \infty}\frac{\C_k}{n_k k^2}=
2\log(b) \,
\hspace{.1in}\mbox{in probability.}
\end{equation}

\noindent{\bf Proof of theorem \ref{cover-bary}} Denote by
$\C_k^+$ the time it takes the walk to cover and return to the
root, and by $R_k$ the number of returns to the root until time
$\C_k^+$. By the remark preceding the proof, $\E\C_k^+ -\E \C_k
\le 4k n_k$, so to prove the theorem it suffices to establish that
\begin{equation} \label{eq:covtree3}
\lim_{k\to \infty}\frac{\E \C_k^+}{n_k k^2}=
2\log(b) \, .
\end{equation}
The expected time to return to the root
is the reciprocal of the root's stationary probability $b/(2n_k-2)$, so by Wald's lemma
\begin{equation} \label{wald}
\E(\C_k^+)=\frac{2n_k-2}{b} \E(R_k) \,.
\end{equation}
Thus the theorem reduces to showing
\begin{equation} \label{eq:R}
\lim_{k\to \infty}\frac{\E\R_k}{k^2} =b\log(b) \, .
\end{equation}

We start by reproducing the straightforward proof of the upper bound.
Denote by
$R_v$ the number of returns to the root of $\Gamma_k$ until the first visit to $v$,
and observe that $R_k$ is the maximum of $R_v$ over all leaves $v$ at level $k$.
At each visit to the root, the chance to hit a specific leaf $v$ before returning to the root is $1/bk$,
whence
$$\P[R_v >rbk^2] \le (1- \frac{1}{bk})^{rbk^2} \le e^{-rk} \,.
 $$
Summing over all leaves, we infer that
\begin{equation} \label{boole}
\P[R_k >rbk^2]   \le \min\{1, b^k e^{-rk}\} \,.
\end{equation}
Integrating over $r>0$,
\begin{equation} \label{boole2}
\E[R_k]   \le  bk^2(\log b + 1/k) \,.
\end{equation}
This yields the upper bound in (\ref{eq:R}).
To prove a lower bound, Aldous~\cite{aldous-t} uses a delicate recursion,
and an embedded branching process argument. Here we will give the shortest argument we know,
which only involves an embedded branching process.
Given $\lambda < \log b$,  our next goal is to show that
\begin{equation} \label{low}
\P[R_k >\lambda bk^2]   \to 1 \mbox{ \rm as } k \to \infty \,.
\end{equation}
Let $T_\lambda$ be the number of steps until the root
is visited $\lambda b k^2$ times.

Fix $r \in (\lambda, \log b)$, and
choose $\ell$ large, depending on $r$.
Let $v$ be a vertex at level $k-(j+1)\ell$ of $\Gamma_k$, and
suppose that $w$ is a descendant of $v$ at level $k-j\ell$.

Observe that the expected number of visits to $v$
by time $T_\lambda$ is $\lambda (b+1) k^2 $, and
the expected number of excursions between
$v$ and $w$ by time $T_\lambda$ is $\lambda k^2 /\ell$.

Say that $w$ is ``special'' if the number of excursions from $v$ to $w$
by time $T_\lambda$ is at most $r\ell j^2$. Note that vertices close to the root
(i.e., at level $k-j\ell$ where $r \ell^2 j^2 > \lambda k^2 $)
are special with high probability, because $r>\lambda$. If $k>(j+2)\ell$, then every visit to $v$
is equally likely to start an excursion to $w$ as to the ancestor of $v$
at distance $\ell$ from $v$. Thus,
if $v$ is special then
$w$ is special with probability at least $\P[X<r \ell j^2]$,
where $X$ has binomial law with parameters $r \ell (j^2+(j+1)^2)$ and $1/2$.
By the central limit theorem, as $j$ grows,
$\P[X<r \ell j^2] \to \P(Z>(2r\ell)^{1/2}]$, where $Z$ is standard normal.
Since $r< \log b$, we find that   $\P(Z>(2r\ell)^{1/2}] >b^{-\ell}$,
if $\ell$ is large enough. Therefore,  special vertices considered at jumps
of $2\ell$ levels (to ensure the required independence) dominate a supercritical branching process;
the survival probability tends to 1 as $k \to \infty$, because vertices near the root
are almost guaranteed to be special.
This establishes (\ref{low}). It follows that $\E(R_k)> \lambda bk^2$ for large $k$,
and since $\lambda <\log b$ is arbitrary, this completes the proof of (\ref{eq:R}) and the theorem.

\noindent{\bf Remark}  The argument above is quite robust: it
readily extends to family trees of Galton watson trees with mean
offspring $b>1$. With a little more work, using the notion of
quasi-Bernoulli percolation (see \cite{lyons} or  \cite{climb}),
it can be extended to the first $k$ levels of any tree $\Gamma$
that has growth and branching number both equal to $b>1$. The
most robust argument, the truncated second moment method used in
\cite{cover}, is too technical to include here.


\section{From trees to Euclidean space} \label{section 3} \setzero
\vskip-5mm \hspace{5mm }

The following ``dictionary'' was offered in
\cite{in} to illustrate the reduction of certain
intersection problems from Euclidean space to trees:

\vskip .2in
\noindent
\begin{tabular}{ll}
{\bf Problem in Euclidean space}&{\bf Corresponding problem on trees} \\
$\bullet$ How many (independent) Brownian&$\bullet$Which branching processes can \\
paths in $\R^d$ can intersect?& have an infinite line of descent?\\
$\bullet$ What is the probability that several&$\bullet$ What is the probability that a \\
random walk paths, started at random&branching process survives \\
 in a  cube of side-length $2^k$, will intersect?& for at least $k$ generations? \\
$\bullet$ Which sets in $\R^3$ contain&$\bullet$ Which trees percolate at \\
double points of Brownian motion?&a fixed threshold $p$? \\
$\bullet$ What is the Hausdorff dimension&$\bullet$ What is the  dimension  \\
of the intersection of a fixed set in& of a percolation cluster \\
$\R^d$ with one or two Brownian paths?&on a general tree?
\end{tabular}

\vskip .2in

The Brownian analogs of Theorems \ref{theo-2lattice} and \ref{ETlattice},
respectively, are given below.
Throughout, denote by $D(x,\eps)$ the disk of radius
$\epsilon$ centered at $x$.

\begin{theorem}[\cite{cover}]
\label{theo-1p}
For Brownian motion $w_{\T}(\cdot)$ in the two-dimensional torus
${\T}^2$, consider the hitting time of a disk,
\[\TT(x,\eps)=\inf \{t>0\,|\,X_t\in  D(x,\eps)\},\]
and the $\ep$-covering time,
\[\CC_\ep=\sup_{x \in \T^2}\TT(x,\eps) \]
which is the amount of time needed for the Wiener sausage of radius
$\ep$ to completely cover ${\T}^2$. Then
\begin{equation}
\lim_{\ep\rar 0}\frac{\CC_\ep}{ \left(\log
\eps\right)^2}=\frac{2}{\pi} \hspace{.6in}\mbox{a.s.\,.}
\label{p.10}
\end{equation}
\end{theorem}

\begin{theorem}[\cite{planar}]  \label{theo-2}
Denote by $\mu_w$ the occupation measure for a planar Brownian motion
$w(\cdot)$ run for unit time.
Then
\begin{equation}
\lim_{\eps\to 0} \sup_{x \in \R^2} \frac{\mu_w (D(x,\eps))}{\eps^2
\left(\log \frac{1}{\eps}\right)^2}=2\,, \hspace{.6in}a.s.\,.
\end{equation}
\end{theorem}
(This was conjectured by
 Perkins and Taylor~\cite{Perkins-Taylor}.)

The basic approach used to prove these results,
 which goes back to Ray, \cite{Ray}, is to control occupation times using
excursions between concentric discs.
The approximate tree structure that is
(implicitly) used arises by considering discs of the same radius $r$
around different centers and varying
$r$; for fixed centers $x,y$, and ``most'' radii $r$ (on a logarithmic
scale) the discs $D(x,r)$ and $D(y,r)$ are either well-separated
(if $r<<|x-y|$) or almost coincide (if $r>>|x-y|$).

\label{lastpage}

\end{document}